\input amstex
\magnification=\magstep1 
\baselineskip=13pt
\documentstyle{amsppt}
\vsize=8.7truein
\CenteredTagsOnSplits \NoRunningHeads
\def\per{\operatorname{per}}
\def\PP{\Cal P}
\def\EE{\bold{E\thinspace}}

\def\RR{\Cal R}
\def\CC{\Cal C}
\def\xx{\bold{x}}
\def\yy{\bold{y}}
\def\sss{\bold{s}}
\def\ttt{\bold{t}}

\def\Pr{\bold{Pr\thinspace }}
\def\HH{\bold{H}}
\topmatter
\title On the number of  matrices and a random matrix  with prescribed row and column sums
and 0-1 entries \endtitle
\author Alexander Barvinok \endauthor
\address Department of Mathematics, University of Michigan, Ann Arbor,
MI 48109-1043, USA \endaddress
\email barvinok$\@$umich.edu \endemail
\thanks This research was partially supported by NSF Grant DMS 0400617 and a
United States - Israel BSF grant 2006377. \endthanks
 \abstract We consider the set $\Sigma(R,C)$ of all $m \times n$ matrices 
 having 0-1 entries and prescribed row sums $R=\left(r_1, \ldots, r_m \right)$ and column sums 
 $C=\left(c_1, \ldots, c_n \right)$. We prove an asymptotic estimate for the cardinality 
 $|\Sigma(R, C)|$ via the solution to a convex optimization problem. We show that 
 if $\Sigma(R, C)$ is sufficiently large, then
 a random matrix $D \in \Sigma(R, C)$ sampled from the uniform probability measure 
 in $\Sigma(R,C)$ with high probability is close to a particular matrix  $Z=Z(R,C)$ that maximizes the sum of entropies
 of entries among all matrices with row sums $R$, column sums $C$ and entries between 0 and 1. 
 Similar
 results are obtained for 0-1 matrices with prescribed row and column sums and 
 assigned zeros in some positions.
 \endabstract
\date November 2009
\enddate
\endtopmatter

\document

\head 1. Introduction and main results \endhead

Matrices with 0-1 entries and prescribed row and column sums is a classical object which 
appears in many branches of pure and applied mathematics. In combinatorics, such matrices 
encode hypergraphs with prescribed degrees of vertices and related structures, see, for 
example, \cite{LW01}. In algebra, certain structural constants in the ring of symmetric functions
and, consequently, in the representation theory of the symmetric and general linear groups
are expressed as numbers of 0-1 matrices with prescribed row and column sums, see 
Chapter 1 of \cite{Ma95}. In statistics, 0-1 matrices with prescribed row and column sums are 
known as {\it binary contingency tables}, see \cite{C+05}. 

Let $R=\left(r_1, \ldots, r_m\right)$ be a positive integer $m$-vector and let 
$C=\left(c_1, \ldots, c_n \right)$ be a positive integer $n$-vector such that 
$$\split &\qquad \qquad \qquad \sum_{i=1}^m r_i =\sum_{j=1}^n c_j =N  \quad \text{and} \\
&0 < r_i < n \quad \text{for} \quad i=1, \ldots, m \quad \text{and} \quad 0< c_j < m \quad 
\text{for} \quad j=1, \ldots, n. \endsplit$$
Let 
$\Sigma(R,C)$ be the set of all $m \times n$ matrices (binary contingency tables)
$D=\left(d_{ij}\right)$ such that 
$$\split \sum_{j=1}^n d_{ij}=r_i \quad \text{for} \quad i=1, \ldots, m, &\quad \sum_{i=1}^m d_{ij} =c_j
\quad \text{for} \quad j=1, \ldots, n \quad \text{and} \\
& d_{ij} \in \{0,1\}. \endsplit$$
In words: $\Sigma(R,C)$ is the set of 0-1 matrices with row sums $R$ and column sums $C$.
Vectors $R$ and $C$ are called {\it margins} of a matrix $D \in \Sigma(R, C)$.

Our first main result provides an estimate of the cardinality of $\Sigma(R, C)$.
\proclaim{(1.1) Theorem} Let us define the function 
$$\split &F(\xx, \yy) =\left( \prod_{i=1}^m x_i^{-r_i} \right) \left( \prod_{j=1}^n y_j^{-c_j} \right) 
\left( \prod_{ij} (1+x_i y_j) \right) \\  &\text{for} \quad \xx=\left(x_1, \ldots, x_n \right) 
\quad \text{and} \quad \yy=\left(y_1, \ldots, y_n \right) \endsplit $$
and let 
$$\alpha(R, C) =\inf \Sb x_1, \ldots, x_m > 0 \\ y_1, \ldots, y_n > 0 \endSb F(\xx, \yy).$$
Then for the number $|\Sigma(R,C)|$ of $m \times n$ zero-one matrices with row sums $R$ 
and column sums $C$ we have 
$$ \alpha(R, C) \ \geq \ |\Sigma(R, C)| \ \geq
 \ {(mn)! \over (mn)^{mn}}\left(\prod_{i=1}^m {(n-r_i)^{n-r_i} \over (n-r_i)!}
\right) \left( \prod_{j=1}^n {c_j^{c_j} \over c_j!} \right) \alpha(R, C).$$
\endproclaim
Let us estimate the ratio between the lower and the upper bounds for $|\Sigma(R,C)|$
using Stirling's formula
$$s! s^{-s}   =e^{-s} \sqrt{2 \pi s} \left(1+O(s^{-1})\right).$$
Since $$e^{-mn} \left( \prod_{i=1}^m e^{n-r_i} \right) \left( \prod_{j=1}^n e^{c_j} \right) =1,$$
the `` $e^{-s}$ '' contributions from Stirling's formula cancel each other out and 
we obtain 
$$ \alpha(R, C) \ \geq \ |\Sigma(R, C)| \ \geq
 \  (mn)^{-\gamma(m+n)} \alpha(R, C)$$
for some absolute constant $\gamma>0$. 

We note that in many interesting cases 
we have  $|\Sigma(R, C)| =2^{\Omega(mn)}$,
see also Section 3.1,  in which case the estimate of Theorem 1.1 captures 
the logarithmic order of $|\Sigma(R, C)|$.

Let us substitute $x_i =e^{s_i}$, $y=e^{t_i}$ in $F(\xx, \yy)$. Then 
$\ln F(\xx, \yy)=G(\sss, \ttt)$, where
$$\split G(&\sss, \ttt)=-\sum_{i=1}^m r_i s_i - \sum_{j=1}^n c_j t_j +\sum_{ij} \ln \left(1+e^{s_i +t_j}\right)
\\ \text{for} \quad &\sss=\left(s_1, \ldots, s_m\right) \quad \text{and} \quad 
\ttt=\left(t_1, \ldots, t_n\right). \endsplit$$
One can observe that $G(\sss, \ttt)$ is a convex function on ${\Bbb R}^m \times {\Bbb R}^n$,
hence to compute the infimum of $G(\sss, \ttt)$ one can use any of the efficient convex 
optimization algorithms, see, for example, \cite{NN94}.

Suppose that margins $R, C$ are such that the set $\Sigma(R, C)$ is not 
empty and let us consider $\Sigma(R, C)$ as a finite probability space with the uniform 
measure. Let us pick a random matrix $D \in \Sigma(R, C)$. What is $D$ likely to look like?
This question is of some interest to statistics: a binary contingency table $D=\left(d_{ij}\right)$
may represent certain statistical data (for example, $d_{ij}$ may be equal to 1 or 0 depending 
on whether or not  Darwin finches of the $i$-th species can be found on the $j$-th Galapagos island, as in \cite{C+05}). One can 
condition on the row and column sums and ask what is special about a particular table 
$D \in \Sigma(R, C)$, considering all tables in $\Sigma(R, C)$ as equiprobable, see \cite{C+05}.
To answer this question we need to know what a random table $D \in \Sigma(R, C)$ looks like.

We prove that with high probability $D$ is close to a particular matrix 
$Z$ with row sums $R$ and column sums $C$ and entries between 0 and 1, 
which we call the {\it maximum entropy matrix}.

\subhead (1.2) The maximum entropy matrix \endsubhead
For $0 \leq x \leq 1$ let us consider the {\it entropy function}
$$H(x)=x \ln {1 \over x} +(1-x) \ln {1 \over 1-x}.$$
As is known, $H$ is  a strictly concave function with $H(0)=H(1)=0$.

For an $m \times n$ matrix $X=\left(x_{ij}\right)$ such that $0 \leq x_{ij} \leq 1$ for all 
$i,j$, we define 
$$H(X) = \sum_{ij} H\left(x_{ij}\right).$$
Assume that $\Sigma(R,C)$ is non-empty.
Let us consider the polytope $\PP(R, C)$ of matrices $X=\left(x_{ij}\right)$
such that
$$\split \sum_{j=1}^n x_{ij}=r_i \quad \text{for} \quad i=1, \ldots, m, & \quad 
\sum_{i=1}^m x_{ij}=c_j \quad \text{for} \quad j=1, \ldots, n \quad \text{and} \\
&0 \leq x_{ij} \leq 1 \quad \text{for all} \quad i,j. \endsplit$$
Since $H(X)$ is strictly concave, it attains a unique maximum $Z=Z(R, C)$ on $\PP(R, C)$,
 which we call the {\it maximum entropy matrix} with margins $(R, C)$.
 
 For example, if all $r_i$ are equal, then by the symmetry argument we must have $Z=\left(z_{ij}\right)$
 where $z_{ij}=c_j/m$ for all $i, j$.
 
 The following observation characterizes the 
 maximum entropy matrix as the solution to the problem that is convex dual to the 
 optimization problem of Theorem 1.1.

\proclaim{(1.3) Lemma} Suppose that the polytope $\PP(R, C)$ has a non-empty interior,
that is, contains a matrix 
$Y=\left(y_{ij}\right)$ such that $ 0 < y_{ij} <1$ for all $i,j$. Then the 
infimum $\alpha(R, C)$ in Theorem 1.1 is attained at a particular
point $\xx^{\ast}=\left(\xi_1, \ldots, \xi_m \right)$ and $\yy^{\ast}=\left(\eta_1, \ldots, \eta_n\right)$.
For the maximum entropy matrix $Z=\left(z_{ij}\right)$ we have
$$z_{ij}={\xi_i \eta_j \over 1+\xi_i \eta_j} \quad \text{for all} \quad i,j \tag1.3.1$$
and, moreover,
$$\alpha(R, C)=e^{H(Z)}. \tag1.3.2$$
Conversely, if the infimum $\alpha(R, C)$ in Theorem 1.1 is attained at a certain point
$\xx^{\ast}=\left(\xi_1, \ldots, \xi_m \right)$ and $\yy^{\ast}=\left(\eta_1, \ldots, \eta_n\right)$
then for the maximum entropy matrix  $Z=\left(z_{ij}\right)$ equations (1.3.1) and (1.3.2) hold.
\endproclaim

The condition that the polytope $\PP(R, C)$ has a non-empty interior 
 is equivalent to the requirement that for every choice of 
$1 \leq k \leq m$ and $1 \leq l \leq n$ there is a matrix $D^0 \in \Sigma(R, C)$, 
$D^0=\left(d_{ij}^0\right)$, such
that $d_{kl}^0=0$ and there is a matrix $D^1 \in \Sigma(R, C)$, $D^1=\left(d_{ij}^1 \right)$,
such that $d_{kl}^1 =1$. One can take $Y$ to be the average of all matrices $D \in \Sigma(R, C)$. 
In other words, we require the set $\Sigma(R, C)$ to be reasonably large.
We also observe that if $r_i c_j < N$ for all $i,j$ (recall that $N$ is the total
sum of the matrix entries) one can choose $y_{ij}=r_i c_j/N$.
\bigskip

We prove that with high probability a random matrix $D \in \Sigma(R, C)$ is close to the maximum entropy matrix 
$Z$ as far as sums over subsets of entries are concerned.

For a subset 
$$S \subset \Bigl\{(i,j): \quad i=1, \ldots, m,\quad  j=1,\ldots, n \Bigr\}$$ and an 
$m \times n$ matrix $A=\left(a_{ij}\right)$, let us denote 
$$\sigma_S(A)=\sum \Sb (i,j) \in S \endSb a_{ij},$$
the sum of the entries of $A$ indexed by $S$.

In what follows, we are interested in the case of the density $N/mn$ separated from 0.
Without loss of generality, we assume that $n \geq m$.

\proclaim{(1.4) Theorem} 
Let us fix numbers $\kappa >0$ and $0 < \delta< 1$.
Then there exists a number $q=q(\kappa, \delta)$ such that the following holds. 

Let $(R, C)$ be margins such that $n \geq m > q$ and the polytope 
$\PP(R, C)$ has a non-empty interior, and let $Z \in \PP(R, C)$ be the maximum entropy matrix.
Let $S \subset \bigl\{(i,j): \quad i=1,\ldots, m,\  j=1, \ldots, n \bigr\}$ be a subset such that 
$\sigma_S(Z) \ \geq \ \delta mn$
and let 
$$\epsilon= \delta{\ln n  \over \sqrt{m}}.$$
If $\epsilon \leq 1$ then 
$$\split \Pr \Bigl\{D \in &\Sigma(R, C): \\ &(1-\epsilon) \sigma_S(Z) \ \leq \ \sigma_S(D) \ \leq \ 
(1+ \epsilon) \sigma_S(Z) \Bigr\}  \geq 1-2n^{-\kappa n}. \endsplit$$
\endproclaim

Let us associate with a non-negative, non-zero $m \times n$ matrix $A=\left(a_{ij}\right)$ a finite 
probability space on the ground set $\{(i,j):\ i=1, \ldots, m, \ j=1, \ldots, n\}$ with 
$\Pr\{(i,j)\}=a_{ij}/N$, where $N >0$ is the total sum of matrix entries. Theorem  1.4 asserts that 
the probability space associated with the maximum entropy matrix $Z$ reasonably well approximates 
the probability space associated with a random binary contingency table $D \in \Sigma(R,C)$ as 
far as events $S$ whose probability is separated from 0 are concerned.

The following interpretation of the maximum entropy matrix was suggested to the author by 
J.A. Hartigan, see \cite{BH09}.

\proclaim{(1.5) Theorem} Let $Z=\left(z_{ij}\right)$ be the $m \times n$ maximum entropy matrix 
with margins $(R, C)$ and let us suppose that the polytope $\PP(R, C)$ has a non-empty interior.
Let $X=\left(x_{ij}\right)$ be the random $m \times n$ matrix of independent Bernoulli 
random variables such that $$\EE X=Z.$$ In other words, $\Pr\left\{x_{ij} =1 \right\}=z_{ij}$ and 
$\Pr\left\{x_{ij}=0\right\}=1-z_{ij}$ independently for all $i, j$. 
Then the probability mass function of $X$ is constant on the set $\Sigma(R, C)$ of binary 
contingency tables with margins $(R, C)$, and, moreover, 
$$\Pr\bigl\{X=D\bigr\}=e^{-H(Z)} \quad \text{for all} \quad D \in \Sigma(R, C).$$
\endproclaim 

The distribution of the random matrix $X$ in Theorem 1.5 can be characterized as the maximum 
entropy distribution in the class consisting of all probability distributions on the set $\{0, 1\}^{m \times n}$
of matrices with 0-1 entries whose expectations lie in the affine subspace consisting of the matrices 
with row sums $R$ and column sums $C$, see \cite{BH09}.

\head 2. Extensions and ramifications \endhead

Our results hold in a somewhat greater generality. Let us fix an $m \times n$ non-negative 
matrix $W=\left(w_{ij}\right)$, which we call the matrix of {\it weights}.
Let us consider the following {\it partition function}
$$|\Sigma(R, C; W)| =\sum \Sb D \in \Sigma(R, C) \\ D=\left(d_{ij}\right) \endSb 
\prod \Sb i,j  \\ d_{ij}=1 \endSb w_{ij}.$$
In particular, if $w_{ij}=1$ for all $i,j$ then 
$|\Sigma(R, C; W)|=|\Sigma(R, C)|$. If $w_{ij}  \in \{0, 1\}$ then the partition 
function counts binary contingency tables with zeros assigned to some positions: 
the value of $|\Sigma(R, C; W)|$ is equal to the number of $m \times n$ matrices 
$D=\left(d_{ij}\right)$ such that the row sums of $D$ are $R$, the column sums of $D$ are $C$,
$d_{ij} \in \{0, 1\}$ for all $i,j$, and, additionally, $d_{ij}=0$ if $w_{ij}=0$. In combinatorial terms,
the set $\Sigma(R, C; W)$ can be interpreted as the set of all subgraphs with prescribed degrees 
of vertices of a given bipartite graph. Binary contingency tables with preassigned zeros are of 
interest in statistics, see \cite{C+05}.

We prove the following result.
\proclaim{(2.1) Theorem} Let us define the function 
$$\split &F(\xx, \yy; W) =\left( \prod_{i=1}^m x_i^{-r_i} \right) \left( \prod_{j=1}^n y_j^{-c_j} \right) 
\left( \prod_{ij} (1+w_{ij} x_i y_j) \right) \\  &\text{for} \quad \xx=\left(x_1, \ldots, x_n \right) 
\quad \text{and} \quad \yy=\left(y_1, \ldots, y_n \right) \endsplit $$
and let 
$$\alpha(R, C; W) =\inf \Sb x_1, \ldots, x_m > 0 \\ y_1, \ldots, y_n > 0 \endSb F(\xx, \yy; W).$$
Then for the partition function $|\Sigma(R,C; W)|$  we have 
$$ \split \alpha(R, C; W) \ \geq \ |\Sigma(R, &C; W)| \\  \geq
 &{(mn)! \over (mn)^{mn}}\left(\prod_{i=1}^m {(n-r_i)^{n-r_i} \over (n-r_i)!}
\right) \left( \prod_{j=1}^n {c_j^{c_j} \over c_j!} \right) \alpha(R, C; W). \endsplit$$
\endproclaim
As before, the function obtained as the result 
of the substitution $x_i=e^{t_i}$, $y_j=e^{s_j}$ in $\ln F(\xx, \yy; W)$,
$$\split G(&\sss, \ttt; W)=-\sum_{i=1}^m r_i s_i - \sum_{j=1}^n c_j t_j +\sum_{ij} \ln \left(1+w_{ij}e^{s_i +t_j}\right)
\\ \text{for} \quad &\sss=\left(s_1, \ldots, s_m\right) \quad \text{and} \quad 
\ttt=\left(t_1, \ldots, t_n\right) \endsplit$$
is convex on ${\Bbb R}^m \times {\Bbb R}^n$, hence computing $\alpha(R, C; W)$ is a convex 
optimization problem.

Let us assume now that $w_{ij} \in \{0, 1\}$ for all $(i,j)$ and let us consider the
set $\Sigma(R, C; W)$ of all $m \times n$ binary contingency tables $D=\left(d_{ij}\right)$ with the additional constraint that $d_{ij}=0$ if $w_{ij}=0$. Assuming that $\Sigma(R, C; W)$ is 
not empty, we consider this set as a finite probability space with the uniform measure. We 
call matrix $W$ the {\it pattern}.
We are interested in what a random table $D \in \Sigma(R, C; W)$ looks like.
We define the maximum entropy matrix as before.
\subhead (2.2) The maximum entropy matrix \endsubhead
Suppose that the set $\Sigma(R,C; W)$ is non-empty.
Let us consider the polytope $\PP(R,C; W)$ of $m \times n$ matrices $X=\left(x_{ij}\right)$
such that
$$\split \sum_{j=1}^n & x_{ij}=r_i \quad \text{for} \quad i=1, \ldots, m, \quad 
\sum_{i=1}^m x_{ij}=c_j \quad \text{for} \quad j=1, \ldots, n,  \\
&0 \leq x_{ij} \leq 1 \quad \text{for all} \quad i,j \quad \text{and} \quad x_{ij}=0 \quad \text{whenever} 
\quad w_{ij}=0.\endsplit$$
Thus $\PP(R, C; W)$ is a face of polytope $\PP(R, C)$ of Section 1.2.
 
 Let $H(X)$ be the entropy function of Section 1.2.
Since $H(X)$ is strictly concave, it attains a unique maximum $Z=Z(R, C; W)$ on 
polytope $\PP(R, C; W)$, which we call 
the {\it maximum entropy matrix} with margins $(R, C)$ and pattern $W$.
\proclaim{(2.3) Lemma} Suppose that the polytope $\PP(R,C;W)$  contains a matrix 
$Y=\left(y_{ij}\right)$ such that $0 < y_{ij}< 1$ whenever $w_{ij}=1$, in which case 
we say that $\PP(R, C; W)$ has a non-empy interior.
Then the infimum $\alpha(R,C;W)$ 
in Theorem 2.1 is attained at a certain point 
$\xx^{\ast}=\left(\xi_1, \ldots, \xi_m \right)$ and  $\yy^{\ast}=\left(\eta_1, \ldots, \eta_n \right)$.
 The maximum entropy matrix
$Z=\left(z_{ij}\right)$ 
satisfies
$$z_{ij}={\xi_i \eta_j \over 1+\xi_i \eta_j} \quad \text{for all} \quad i,j \quad \text{such that} \quad w_{ij}=1
\tag2.3.1$$
Moreover, 
$$\alpha(R, C; W)=e^{H(Z)}. \tag2.3.2$$
Conversely, if the infimum $\alpha(R, C; W)$ is attained at a point $\xx^{\ast}=\left(\xi_1, \ldots, \xi_m \right)$ and $\yy^{\ast}=\left(\eta_1, \ldots, \eta_n \right)$, then for the maximum entropy matrix 
$Z=\left(z_{ij}\right)$ the equations (2.3.1) and (2.3.2) hold.

\endproclaim
For $\PP(R, C; W)$ to have a non-empty interior is equivalent to the requirement that for every 
pair $k,l$ such that 
$w_{kl}=1$ there is a matrix $D^0 \in \Sigma(R, C; W)$, $D^0=\left(d_{ij}^0 \right)$, such that 
$d_{kl}^0=0$ and there is a matrix $D^1 \in \Sigma(R, C; W)$, $D^1=\left(d_{ij}^1 \right)$, 
such that $d_{kl}^1 =1$. 
In other words, we require the set $\Sigma(R, C; W)$ to be reasonably large.

We prove an analogue of Theorem 1.4. We consider subsets
$$S \subset \Bigl\{ (i,j): \quad w_{ij}=1 \Bigr\}.$$
As before, we denote by $\sigma_S(A)$ the sum of the entries of a matrix $A$ indexed by the 
subset $S$. 

\proclaim{(2.4) Theorem} 
Let us fix numbers $\kappa>0$ and $0 < \delta< 1$.
Then there exists a number $q=q(\kappa, \delta)$ such that the following holds. 

Let $(R, C)$ be margins such that $n \geq m > q$ and the polytope 
$\PP(R, C; W)$ has a non-empty interior, and let $Z \in \PP(R, C; W)$ be the maximum entropy matrix.
Let $S \subset \bigl\{(i,j): \quad w_{ij}=1\bigr\}$ be a subset such that 
$\sigma_S(Z) \ \geq \ \delta mn$
 and let 
$$\epsilon= \delta{\ln n \over \sqrt{m}}.$$
If $\epsilon \leq 1$ then 
$$\split \Pr \Bigl\{D \in &\Sigma(R, C; W): \\ &(1-\epsilon) \sigma_S(Z) \ \leq \ \sigma_S(D) 
 \ \leq \ (1+ \epsilon) \sigma_S(Z) \Bigr\}  
\geq 1-2n^{-\kappa n}. \endsplit$$
\endproclaim
The statement of the theorem is, of course, vacuous unless pattern $W$ contains 
$\Omega(mn)$ ones.

There is an analogue of Theorem 1.5.

\proclaim{(2.5) Theorem} Suppose that the polytope $\PP(R, C; W)$ has a non-empty interior 
and let $Z \in \PP(R, C; W)$ be the maximum entropy matrix. 
Let $X=\left(x_{ij}\right)$ be the random $m \times n$ matrix of independent Bernoulli 
random variables such that $$\EE X=Z,$$ that is, $\Pr\left\{x_{ij} =1 \right\}=z_{ij}$ and 
$\Pr\left\{x_{ij}=0\right\}=1-z_{ij}$ independently for all $i, j$. 
Then the probability mass function of $X$ is constant on the set $\Sigma(R, C; W)$ and, moreover, 
$$\Pr\bigl\{X=D\bigr\}=e^{-H(Z)} \quad \text{for all} \quad D \in \Sigma(R, C; W).$$
\endproclaim

\head 3. Comparisons with the literature \endhead

There is a vast literature on  0-1 matrices with prescribed row and column sums and 
with or without zeros in prescribed positions, see for example, Chapter 16 of \cite{LW01}, 
\cite{Ne69}, \cite{Be74}, \cite{GC77},
 recent \cite{CM05}, \cite{G+06}, \cite{C+08}, \cite{GM09} and references 
therein. A simple and efficient criterion for the existence of a 0-1 matrix with 
prescribed row and column sums is given by the classical Gale-Ryser  Theorem; in the 
case of enforced zeros, the question reduces to the existence of a network flow, see 
for example, Chapter 16 of \cite{LW01}. Estimating the number of such matrices also attracted 
a lot of attention.
Precise asymptotic formulas for 
the number of matrices were obtained in sparse cases for which 
$r_i \ll n$ and $c_j \ll m$ \cite{Ne69}, \cite{Be74}, \cite{G+06}, 
the regular case of all row sums $r_i$ equal and all column sums $c_j$ equal \cite{C+08} 
and cases close to regular \cite{C+08}, \cite{GM09}. Formulas of Theorems 1.1 and  2.1
are not so precise but they are applicable to a wide class of margins $(R, C)$ and 
they uncover some interesting features of the numbers $|\Sigma(R, C)|$ and $|\Sigma(R, C; W)|$.

The following construction
provides some insight into the combinatorial interpretation of the number $\alpha(R, C)$ from
Theorem 1.1.
\subhead (3.1) Cloning the margins \endsubhead Let us fix some margins $R, C$ for 
which the set $\Sigma(R, C)$ is not empty, and, moreover, the polytope 
$\PP(R, C)$ contains an interior point, so the conditions of Lemma 1.3 are 
satisfied. Let $R=\left(r_1, \ldots, r_m \right)$ and 
$C=\left(c_1, \ldots, c_n \right)$.
For a positive integer $k$, let us define the $km$-vector 
$$R_k=\left(\underbrace{kr_1, \ldots, kr_1}_{\text{$k$ times}}, \ldots, 
\underbrace{kr_m, \ldots, kr_m}_{\text{$k$ times}}\right)$$
and the $kn$-vector 
$$C_k=\left(\underbrace{kc_1, \ldots, kc_1}_{\text{$k$ times}}, \ldots, 
\underbrace{kc_n, \ldots, kc_n}_{\text{$k$ times}}\right).$$
In other words, we obtain margins $\left(R_k, C_k \right)$ if we choose a matrix 
$Y \in \PP(R, C)$ and then create a new block matrix $Y_k$ by 
arranging $k^2$ copies of $Y$ into a $km \times kn$ matrix.
Then $R_k$ is the vector of row sums of $Y_k$ and $C_k$ is the vector of column sums of $Y_k$.
Clearly, the conditions of Lemma 1.3 are satisfied for $(R_k, C_k)$.

Theorem 1.1 then implies that 
$$\lim_{k \longrightarrow +\infty} \left|\Sigma(R_k, C_k)\right|^{1/k^2}=\alpha(R, C). \tag3.1.1$$
Indeed, the infimum $\alpha(R, C)$ is attained at a certain point
$$\xx^{\ast}=\left(\xi_1, \ldots, \xi_m \right) \quad \text{and} \quad 
\yy^{\ast}=\left(\eta_1, \ldots, \eta_n \right).$$
It is not hard to see that the infimum $\alpha\left(R_k, C_k \right)$ is attained at 
$$\split &\xx^{\ast}_k=\left(\underbrace{\xi_1, \ldots, \xi_1}_{\text{$k$ times}}, \ldots, 
\underbrace{\xi_m, \ldots, \xi_m}_{\text{$k$ times}} \right) \quad \text{and} \\
&\yy^{\ast}_k=\left(\underbrace{\eta_1, \ldots, \eta_1}_{\text{$k$ times}}, \ldots, 
\underbrace{\eta_n, \ldots, \eta_n}_{\text{$k$ times}} \right). \endsplit$$ 

\subhead (3.2) Asymptotic repulsion in the space of matrices \endsubhead
A natural candidate for an approximation of $|\Sigma(R, C)|$ is the ``independence estimate''
$$I(R, C) = {mn \choose N}^{-1} \prod_{i=1}^m {n \choose r_i}  
\prod_{j=1}^n {m \choose c_j}, \tag3.2.1$$
see \cite{GC77}, \cite{G+06}, and \cite{C+08}. 

The intuitive meaning of (3.2.1) is as follows.
Let us consider the set of all $m \times n$ matrices with 0-1 entries and with the total sum of entries 
equal to $N$ as a finite probability space with the uniform measure. Let us consider 
the two events in this space: the event $\RR$ consisting of the matrices with row sums $R$ 
and the event $\CC$ consisting of the matrices with column sums $C$. One can see that 
$$ \Pr (\RR) ={mn \choose N}^{-1} \prod_{i=1}^m {n \choose r_i} \quad \text{and} 
\quad \Pr(\CC)={mn \choose N}^{-1} \prod_{j=1}^n {m \choose c_j}$$
and that 
$$|\Sigma(R, C)| ={mn \choose N} \Pr(\RR \cap \CC).$$
Thus the value of (3.2.1) equals $|\Sigma(R, C)|$ if the events $\RR$ and 
$\CC$ are independent. It turns out that (3.2.1) indeed approximates $|\Sigma(R, C)|$ 
reasonably well in the sparse and near-unform cases, see \cite{G+06} and \cite{C+08}.

However, for generic $R$ and $C$, the independence estimate $I(R, C)$ overestimates 
$|\Sigma(R, C)|$ by a $2^{\Omega(mn)}$ factor. To see why, let us fix some margins
 $R=\left(r_1, \ldots, r_m \right)$ and $C=\left(c_1, \ldots, c_n \right)$
such that not all row sums $r_i$ are equal and not all column sums $c_j$ are equal and 
the conditions of Lemma 1.3 are satisfied.
Let us consider the cloned margins $R_k$ and $C_k$ as in Section 3.1.

Applying Stirling's formula, we get 
$$\split &\lim_{k \longrightarrow +\infty} I\left(R_k, C_k\right)^{1/k^2}\\ 
&\qquad = \exp\biggl\{-mn H\left({N \over mn} \right) +  n \sum_{i=1}^m H\left({r_i \over n} \right) \\ 
&\qquad \qquad \qquad \qquad + m \sum_{j=1}^n H \left({c_j \over m}\right) \biggr\}, \endsplit \tag3.2.2 $$
where 
$H$ is the entropy function, see Section 1.2. To compare (3.2.2) and (3.1.1) we use Lemma 1.3 and the
multivariate entropy function 
$$\HH(p_1, \ldots, p_k) =\sum_{i=1}^k p_k \ln {1 \over p_k},$$
where $p_1, \ldots, p_k$ are non-negative numbers such that $p_1 + \ldots + p_k=1$.
Thus $H(x)=\HH(x, 1-x)$ for $0 \leq x \leq 1$ and we rewrite (3.2.2) as 
$$\split \lim_{k \longrightarrow +\infty} {1 \over k^2} & \ln I\left(R_k, C_k \right) \\
=&N \HH\left({r_1 \over N}, \ldots, {r_m \over N} \right) 
+(mn-N) \HH\left({n-r_1 \over mn- N}, \ldots, {n-r_m \over mn - N} \right) \\
+& N \HH\left({c_1 \over N}, \ldots, {c_n \over N} \right) +(mn-N) 
\HH\left({m-c_1 \over mn- N}, \ldots, {m-c_n \over mn - N} \right) \\ &\quad -N \ln N - (mn-N) \ln (mn-N). \endsplit $$
On the other hand, applying Lemma 1.3, we can rewrite (3.1.1) as 
$$\split \lim_{k \longrightarrow +\infty} {1 \over k^2} &\ln \left| \Sigma\left(R_k, C_k \right) \right|\\
 =&N \HH\left({z_{ij} \over N} \right) + (mn-N) \HH \left(1-z_{ij} \over mn- N\right) - N \ln N \\
&\qquad -(mn-N) \ln (mn-N), \endsplit$$
where $Z=\left(z_{ij}\right)$ is the maximum entropy matrix for margins $(R, C)$.

We now use some classical entropy inequalities, see, for example, \cite{Kh57}.
Namely, by the inequality relating the entropies of two partitions of a probability space and 
the entropy of their intersection, we have 
$$\HH \left({z_{ij} \over N} \right) \leq \HH\left({r_1 \over N}, \ldots, {r_m \over N} \right)  +  
\HH\left({c_1 \over N}, \ldots, {c_n \over N} \right) $$
with the equality if and only if 
$$z_{ij} ={r_i c_j \over N} \quad \text{for all} \quad i,j \tag3.2.3$$
and 
$$\HH \left(1-z_{ij} \over mn- N\right) \leq
  \HH\left({n-r_1 \over mn- N}, \ldots, {n-r_m \over mn - N} \right) + \HH\left({m-c_1 \over mn- N}, \ldots, {m-c_n \over mn - N} \right)$$
  with the equality if and only if 
  $$1-z_{ij} ={(n-r_i) (m -c_j) \over mn - N} \quad \text{for all} \quad i,j. \tag3.2.4$$
  However, if we have both (3.2.3) and (3.2.4), we must have 
  $(r_i m -N)(c_j n- N)=0$, so unless all row sums $r_i$ are equal or all column sums $c_j$ 
  are equal, we have 
  $$\lim_{k \longrightarrow +\infty} \left| \Sigma(R_k, C_k)\right|^{1/k^2}  \quad 
  < \quad \lim_{k \longrightarrow +\infty} I\left(R_k, C_k \right)^{1/k^2}.$$
  Therefore, as $k$ grows, the independence estimate (3.2.1) overestimates the number of 0-1 matrices
  with row sums $R_k$ and column sums $C_k$ by a factor of $2^{\Omega(k^2)}$.
  In probabilistic terms, as $k$ grows, the event $\RR_k$ consisting of the 0-1 matrices with 
  row sums $R_k$ and the event $\CC_k$ consisting of the 0-1 matrices with column sums $C_k$ repel 
  each other (the events are negatively correlated) instead of being asymptotically independent. 
  
  The procedure of cloning described in Section 3.1 produces margins of increasing 
  size with the following features: the density remains separated from 0 and 1, and if the margins 
  were non-uniform initially,  they stay away from uniform. One can show that for more general sequences of margins that 
  share these two features, we have the asymptotic repulsion of the event consisting of 
  the 0-1 matrices with prescribed row sums and the event consisting of the 0-1 matrices
  with prescribed column sums.  This is in contrast to the case of contingency tables (non-negative 
  integer matrices with prescribed row and column sums), where we have the asymptotic 
  attraction of the events \cite{Ba09}.

\subhead  (3.3) Randomized counting and sampling \endsubhead
Jerrum, Sinclair, and Vigoda \cite{J+04} showed how to apply their algorithm for computing the
permanent of a non-negative matrix to construct a fully polynomial randomized approximation
scheme (FPRAS) to compute $|\Sigma(R, C)|$ and, more generally, $|\Sigma(R, C; W)|$,
where $W$ is a 0-1 pattern, see also \cite{B+07}.
 Furthermore, they obtained a polynomial time algorithm for 
sampling a random $D \in \Sigma(R, C)$ and $D \in \Sigma(R, C; W)$ from a ``nearly uniform''
distribution. This problem arises naturally in statistics, see, for example, \cite{C+05}.
The estimates of Theorem 1.1 and Theorem 2.1 are not nearly as precise as those
of \cite{J+04}, but they are deterministic, easily computable, and amenable 
to analysis. Similarly, we do not provide a sampling algorithm but show instead in Theorems 
1.4 and 2.4 what a random matrix is likely to look like.

\subhead (3.4) An open question \endsubhead Theorem 1.5 allows us to interpret Theorem 1.4 as a 
law of large numbers for binary contingency tables: with respect to sums $\sigma_S(D)$ 
for sufficiently ``heavy'' sets $S$ of indices, a random binary contingency table $D \in \Sigma(R, C)$ 
behaves approximately as the matrix of independent Bernoulli random variables whose expectation 
is the maximum entropy matrix $Z=\left(z_{ij}\right)$. Similar concentration results can be obtained 
for other well-behaved functions on binary contingency tables. One can ask whether the distribution of
a {\it particular entry} $d_{ij}$ of a random table $D \in \Sigma(R, C)$ converges in distribution to the 
Bernoulli distribution with expectation $z_{ij}$ as the dimensions $m$ and $n$ of the table 
grow in some regular way, for example, when the margins are cloned as in Section 3.1.
\bigskip
Our approach, based on estimating 
combinatorial quantities via solutions to optimization problems, reminds one of that of Gurvits 
\cite{Gu08}. The appearance of entropy in combinatorial counting problems 
reminds one of recent papers of Cuckler and Kahn \cite{CK09a}, \cite{CK09b}, although methods
and results seem to be quite different.

In the rest of the paper, we prove the results stated in Sections 1 and 2.

\head 4. Preliminaries: permanents and scaling \endhead
   
Let $A=\left(a_{ij}\right)$ be an $n \times n$ matrix. The {\it permanent} of $A$
is defined by the expression
$$\per A=\sum_{\sigma \in S_n} \prod_{i=1}^n a_{i \sigma(i)},$$
where $S_n$ is the symmetric group of all permutations of the set $\{1, \ldots, n\}$.   
The relevance of permanents to us is that both values of $|\Sigma(R, C)|$ and
$|\Sigma(R, C; W)|$ can be expressed as permanents of $mn \times mn$ matrices.
This result is not new, for $|\Sigma(R, C)|$ it was observed, for example, in \cite{JS90}.
For $|\Sigma(R, C; W)|$, where $W$ is a 0-1 pattern, a construction is presented in \cite{J+04}.
We give a general construction for $|\Sigma(R, C; W)|$, where $W$ is an arbitrary matrix,
which is slightly different from that of \cite{J+04}.

\proclaim{(4.1) Lemma} Let us choose margins $R=\left(r_1, \ldots, r_m \right)$, 
$C=\left(c_1, \ldots, c_n \right)$ and an $m \times n$ matrix $W=\left(w_{ij}\right)$ of weights.
Let us construct an $mn \times mn$ matrix $A=A(R, C; W)$ as follows.

The rows of $A$ are split into disjoint
\smallskip
$m$ blocks having $n-r_1, \ldots, n-r_m$ rows respectively 
(blocks of type I) 

\noindent and 

$n$ blocks having  $c_1, \ldots, c_n$ rows respectively (blocks of type II).
\smallskip
The columns of $A$ are split into $m$ disjoint blocks of  $n$ columns in each. 
\smallskip
For $i=1, \ldots, m$ the entry of $A$ that lies in a row from the $i$-th block of rows of type I and in a 
column from the $i$-th block of columns is equal to 1.

For $i=1, \ldots, m$ and $j=1, \ldots, n$ the entry of $A$ that lies in a row from the $j$-th block 
of rows of type II and the $j$-th column from the $i$-th block of columns is equal to $w_{ij}$.

All other entries of $A$ are 0s.

Then 
$$|\Sigma(R, C; W)| =\left( \prod_{i=1}^m {1 \over (n-r_i)!} \right) \left(\prod_{j=1}^n {1 \over c_j!} \right)
\per A.$$
\endproclaim
\demo{Proof}
First, we express $|\Sigma(R,C;W)|$ as a coefficient in a certain polynomial.
Let $x_1, \ldots, x_n$ be formal variables and let 
$$e_r\left(x_1, \ldots, x_n \right)=\sum \Sb 1 \leq i_1 < \ldots < i_r \leq n \endSb x_{i_1} \cdots x_{i_r}$$
be the elementary symmetric polynomial of degree $r$. Thus  $e_r\left(x_1, \ldots, x_n \right)$ is 
$$\text{the coefficient of} \quad t^{n-r} \quad \text{in the product} \quad \prod_{j=1}^n \left(t + x_j\right).$$
We observe that $|\Sigma(R,C; W)|$ is
$$\text{the coefficient of} \quad \prod_{j=1}^n x_j^{c_j} \quad \text{in the product} \quad \prod_{i=1}^m e_{r_i} \left(w_{i1} x_1, \ldots, w_{in} x_n \right).$$ 

Summarizing, we conclude that $|\Sigma(R, C; W)|$ is 
$$\text{the coefficient of} \quad \prod_{i=1}^m t_i^{n-r_i} \prod_{j=1}^n x_j^{c_j} 
\quad \text{in the product} \quad \prod_{i=1}^m \prod_{j=1}^n \left(t_i + w_{ij} x_j \right).$$
To express the last coefficient as the permanent of a matrix, we use a convenient scalar product
in the space of polynomials, see, for example, \cite{Ba96} and \cite{Ba07}. Namely, for monomials
$$\xx^{a} =x_1^{\alpha_1} \cdots x_n^{\alpha_n} \quad \text{where} \quad 
a=\left(\alpha_1, \ldots, \alpha_n \right) \quad \text{and} \quad \xx=\left(x_1, \ldots, x_n\right)$$
 we define
$$\langle \xx^a , \xx^b \rangle = \cases \alpha_1! \cdots \alpha_n! &\text{if \ } 
a=b=\left(\alpha_1, \ldots, \alpha_n \right) \\ 0 &\text{if\ } a \ne b \endcases$$
and then extend the scalar product $\langle \cdot, \cdot \rangle$ by bilinearity. Equivalently, 
the scalar product can be defined as follows: let us identify 
${\Bbb R}^n \oplus {\Bbb R}^n={\Bbb C}^n$ via $x+iy=z$ and let $\nu_n$ be the Gaussian measure 
on ${\Bbb C}^n$ with the density 
$$\pi^{-n} e^{-\|z\|^2} \quad \text{where} \quad \|z\|^2=\|x\|^2 + \|y\|^2 \quad \text{for} \quad z=x+iy.$$
Then, for polynomials $f$ and $g$ we have
$$\langle f, g \rangle =\int_{{\Bbb C}^n} f(z) \overline{g(z)} \ d\nu_n,$$
where $\overline{g}$ is the complex conjugate of $g$, see for example, Section 4 of \cite{Ba07}.

The convenient 
property of the scalar product is that if 
$$p(\xx)=\prod_{l=1}^m \sum_{k=1}^n b_{lk} x_k \quad \text{and} \quad 
q(\xx)=\prod_{l=1}^m \sum_{k=1}^n c_{lk} x_k$$
are products of linear forms, then
$$\langle p, q \rangle = \per D,$$
where $D=\left(d_{ij}\right)$ is the $m \times m$ matrix defined by 
$$d_{ij}=\sum_{k=1}^n b_{ik} c_{j k} \quad \text{for all} \quad i,j,$$
see Lemma 4.5 of \cite{Ba07} or, for a more general identity, Theorem 3.8 of \cite{Gu04}.
Thus we may write
$$\split |\Sigma(R, C; W)|=& \left(\prod_{i=1}^m {1 \over (n-r_i)!}\right) \left( \prod_{j=1}^n {1 \over c_j!} \right) \\
&\qquad \qquad \times \left\langle \prod_{i=1}^m t_i^{n-r_i} \prod_{j=1}^n x_j^{c_j}, \quad \prod_{i=1}^m \prod_{j=1}^n
\left(t_i + w_{ij} x_j \right) \right\rangle \\
= & \left(\prod_{i=1}^m {1 \over (n-r_i)!}\right) \left( \prod_{j=1}^n {1 \over c_j!}  \right)
\per A.\endsplit $$
 {\hfill \hfill \hfill} \qed
\enddemo

\subhead (4.2) Matrix scaling and the van der Waerden bound \endsubhead
Let $B=\left(b_{ij}\right)$ be an $n \times n$ matrix. Matrix $B$ is called {\it doubly stochastic}
if
$$\split \sum_{j=1}^n &b_{ij}=1 \quad \text{for} \quad i=1, \ldots, m, \quad \sum_{i=1}^n b_{ij} =1 
\quad \text{for} \quad j=1, \ldots n \quad \text{and} \\
&b_{ij} \geq 0 \quad \text{for all} \quad \text{for all} \quad i,j. \endsplit$$
The classical bound conjectured by van der Waerden and proved by Falikman and Egorychev,
see Chapter 12 of \cite{LW01} and also \cite{Gu08} for exciting new developments, states that 
$$\per B \geq {n! \over n^n}$$
if $B$ is a doubly stochastic matrix.

Linial, Samorodnitsky, and Wigderson \cite{L+00} introduced the following very useful 
{\it scaling} method of 
approximating permanents of non-negative matrices. Given a non-negative $n \times n$ matrix 
$A=\left(a_{ij}\right)$ 
one finds non-negative numbers $\lambda_1, \ldots, \lambda_n$ and $\mu_1, \ldots, \mu_n$ 
and a doubly stochastic matrix $B=\left(b_{ij}\right)$ such that 
$$a_{ij}=\lambda_i \mu_j b_{ij} \quad \text{for all} \quad i,j.$$
Then 
$$\per A = \left( \prod_{i=1}^n \lambda_i \right) \left( \prod_{j=1}^n \mu_j \right) \per B$$
and an estimate of $\per B$ (such as the van der Waerden estimate) implies an estimate of 
$\per A$. If $A$ is strictly positive, such doubly stochastic matrix $B$ and scaling factors 
$\lambda_i$, $\mu_j$ always exist. In our situation, matrix $A$  constructed in Lemma 4.1 
is only non-negative. We will not always be able to scale it to a doubly stochastic matrix $B$ exactly, 
but we will scale it {\it approximately}.

We restate a weaker form of Proposition 5.1 from \cite{L+00} regarding almost doubly stochastic 
matrices.

\proclaim{(4.3) Lemma} For any $n$ there exists an $\epsilon_0 =\epsilon_0(n) >0$ and 
a function $\phi(\epsilon)$, $0 < \epsilon < \epsilon_0$, such that 
$$\lim_{\epsilon \longrightarrow 0+} \phi(\epsilon) =1$$ 
and for any $n \times n$ non-negative matrix $B=\left(b_{ij}\right)$ such that 
$$\sum_{i=1}^n b_{ij} =1 \quad \text{for} \quad j=1, \ldots, n$$
and 
$$1-\epsilon \ \leq \ \sum_{j=1}^n b_{ij} \ \leq \ 1+ \epsilon \quad \text{for} \quad i=1, \ldots, n$$
for some $0 \leq \epsilon < \epsilon_0$,
we have 
$$\per B \geq {n! \over n^n} \phi(\epsilon).$$
\endproclaim 
{\hfill \hfill \hfill} \qed

From \cite{L+00}, one can choose $\epsilon_0=1/n$ and $\phi(\epsilon)=(1-\epsilon n)^n$.

\head 5. Proofs of Theorems 1.1 and 2.1 \endhead

We prove Theorem 2.1 only since Theorem 1.1 is a particular case of Theorem 2.1.
We start with a straightforward observation.

\proclaim{(5.1) Lemma} We have 
$$\split \prod_{ij} \left(1+ w_{ij} x_i y_j \right) =\sum_{(R, C)} |\Sigma(R, &C; W)|  \xx^R \yy^C,
\quad \text{where} \\ &\xx^R =x_1^{r_1} \cdots x_m^{r_m}, \quad 
\yy^C=y_1^{c_1} \cdots y_n^{c_n}, \endsplit$$
and the sum is taken over all margins $R, C$.
\endproclaim
{\hfill \hfill \hfill} \qed

Next, we need a technical lemma.
\proclaim{(5.2) Lemma} Let $W=\left(w_{ij}\right)$ be an $m \times n$ non-negative matrix 
such that 
$$\alpha(R, C; W) >0.$$
Then, for any $\epsilon > 0$ there exist points $\xx=\xx(\epsilon)$ and $\yy=\yy(\epsilon)$,
$\xx=\left(x_1, \ldots, x_n \right)$ and $\yy=\left(y_1, \ldots, y_n \right)$, such that 
$$\split &\left|-r_i + \sum_{j=1}^n {w_{ij} x_i y_j \over 1+ w_{ij} x_i y_j} \right| < \epsilon \quad \text{for} 
\quad i=1, \ldots, m \\
&\left| -c_j + \sum_{i=1}^m {w_{ij} x_i y_j \over 1 + w_{ij} x_i y_j} \right|< \epsilon \quad 
\text{for} \quad j=1, \ldots, n \quad \text{and} \\
& \qquad \qquad x_i, y_j > 0 \quad \text{for all} \quad i,j. \endsplit$$
\endproclaim

\demo{Proof}
Let us consider the function
$$\split G(&\sss, \ttt; W)=-\sum_{i=1}^m r_i s_i - \sum_{j=1}^n c_j t_j +\sum_{ij} \ln \left(1+w_{ij}e^{s_i +t_j}\right)
\\ \text{for} \quad &\sss=\left(s_1, \ldots, s_m\right) \quad \text{and} \quad 
\ttt=\left(t_1, \ldots, t_n\right). \endsplit$$
Then $G(\sss, \ttt; W)$ is convex and 
$$\inf \Sb \sss \in {\Bbb R}^m \\ \ttt \in {\Bbb R}^n \endSb G(\sss, \ttt) = \ln \alpha(R, C; W)  >
-\infty.$$
Hence $G(\sss, \ttt)$ is bounded from below, it is also easy to check that the 
Hessian of $G$ remains bounded on ${\Bbb R}^m \times {\Bbb R}^n$. 
Therefore, the gradient of $G(\sss, \ttt)$  can get arbitrarily close to 0.
That is, for any $\epsilon >0$ there are points
$$\sss(\epsilon)=\left(s_1(\epsilon), \ldots, s_m(\epsilon)\right) \quad \text{and} \quad 
\ttt(\epsilon)=\left(t_1(\epsilon), \ldots, t_n(\epsilon) \right)$$
such that 
$$\split &\left| {\partial \over \partial s_i} G(\sss, \ttt) \big|_{\sss=\sss(\epsilon), \ttt=\ttt(\epsilon)} \right| <
\epsilon \quad \text{for} \quad i=1, \ldots, m \quad \text{and} \\
 &\left| {\partial \over \partial t_j} G(\sss, \ttt) \big|_{\sss=\sss(\epsilon), \ttt=\ttt(\epsilon)} \right| <
\epsilon \quad \text{for} \quad j=1, \ldots, n \endsplit$$
(it suffices to choose $\sss(\epsilon)$ and $\ttt(\epsilon)$ so that the value of 
$G(\sss(\epsilon), \ttt(\epsilon))$ is sufficiently close to the infimum).
In other words,
$$\split &\left| -r_i + \sum_{j=1}^n 
{w_{ij} e^{s_i(\epsilon)+t_j(\epsilon)}  \over 1+ w_{ij} e^{s_i(\epsilon) +t_j(\epsilon)}} \right| 
< \epsilon \quad \text{for} \quad i=1, \ldots, m  \\ &\qquad \qquad \text{and} \\
&\left| -c_j + \sum_{i=1}^m
{w_{ij} e^{s_i(\epsilon) + t_j(\epsilon)} \over 1+ w_{ij} e^{s_i(\epsilon) + t_j(\epsilon)}} \right| 
< \epsilon \quad \text{for} \quad j=1, \ldots, n. \endsplit 
$$
We now let 
$$\split &x_i=x_i(\epsilon)=e^{s_i(\epsilon)} \quad \text{for} \quad i=1, \ldots, m \quad \text{and} \\
&y_j=y_j(\epsilon) =e^{t_j(\epsilon)} \quad \text{for} \quad j=1, \ldots, n. \endsplit$$
{\hfill \hfill \hfill} \qed
\enddemo

\subhead (5.3) Proof of Theorem 2.1 \endsubhead
The upper bound 
$$\alpha(R, C; W) \geq |\Sigma(R, C; W)|$$
follows from Lemma 5.1. Let us prove the lower bound.

If $\alpha(R,C; W)=0$ then $|\Sigma(R,C;W)|=0$ and the lower bound follows.
Hence we assume that $\alpha(R,C;W)>0$. 

Let $A=A(R, C; W)$ be the $mn \times mn$ block matrix constructed in Lemma 4.1.
Let us consider the $mn \times mn$ block matrix $B(\epsilon)$ obtained from $A$ as follows. 
For $\epsilon >0$, let $\xx(\epsilon)=\left(x_1, \ldots, x_m \right)$ and 
$\yy(\epsilon)=\left(y_1, \ldots, y_n \right)$ be the point constructed in Lemma 5.2.
\smallskip
For $i=1, \ldots, m$ we multiply every row of $A$ in the $i$-th block of type I by 
$${1 \over x_i (n-r_i)};$$

For $j=1, \ldots, n$ we multiply every row of $A$ in the $j$-th  block of type II by 
$${y_j \over c_j} \quad \text{for} \quad j=1, \ldots, n;$$

For $i=1, \ldots, m$ and $j=1, \ldots, n$ we multiply the $j$-th column in the $i$-th block of columns
of $A$ by 
$${x_i \over 1+w_{ij} x_i y_j}.$$
This choice of scaling factors is, basically, a lucky guess made in the hope to match the structure
of the function $F(\xx, \yy; W)$.
\smallskip

Thus we have 
$$\split \per A =\left( \prod_{i=1}^m x_i^{n-r_i} (n-r_i)^{n-r_i} \right) &
\left( \prod_{j=1}^n y_j^{-c_j} c_j^{c_j} \right) \\ \times 
&\left( \prod_{ij}  x_i^{-1} \left(1+w_{ij} x_i y_j \right)\right) \per B(\epsilon) \endsplit$$
and hence
$$\split |\Sigma(R,C; W)| =&\left( \prod_{i=1}^m {(n-r_i)^{n-r_i} \over (n-r_i)!} \right) 
\left( \prod_{j=1}^n {c_j^{c_j} \over c_j^{c_j}} \right) \\ &\qquad \qquad \times F\bigl(\xx(\epsilon), \yy(\epsilon); W\bigr) \per B(\epsilon) \\ \geq & \left( \prod_{i=1}^m {(n-r_i)^{n-r_i} \over (n-r_i)!} \right) \left( \prod_{j=1}^n {c_j^{c_j} \over c_j^{c_j}} \right) \\
&\qquad \qquad \times  \alpha(R, C; W) \per B(\epsilon). \endsplit \tag5.3.1$$

Finally, we claim that $B(\epsilon)$ is close to a doubly stochastic matrix.
Indeed,
\smallskip
For $i=1, \ldots, m$ and $j=1, \ldots, n$ the entry of $B(\epsilon)$ that lies in a row from the $i$-th block of rows of type I and in the $j$-th column  from the $i$-th block of columns is equal to 
$${1 \over (n-r_i) (1+w_{ij} x_i y_j)}.$$

For $i=1, \ldots, m$ and $j=1, \ldots, n$ the entry of $B(\epsilon)$ that lies in a row from the $j$-th block 
of rows of type II and the $j$-th column from the $i$-th block of columns is equal to 
$${w_{ij} x_i y_j \over c_j(1+w_{ij} x_i y_j)}.$$

All other entries of $B(\epsilon)$ are 0s.
\smallskip
Let us compute the row sums of $B(\epsilon)$. 

For a row in the $i$-th block of rows of type I the sum equals 
$$ a_i=\sum_{j=1}^n {1 \over (n-r_i)(1+w_{ij} x_i y_j)}.$$
Since 
$$\sum_{j=1}^n {1 \over 1+w_{ij} x_i y_j} =\sum_{j=1}^n {1+ w_{ij} x_i y_j \over 1+ w_{ij} x_i y_j} 
-\sum_{j=1}^n {w_{ij} x_i y_j \over 1+ w_{ij} x_i y_j},$$
by the inequalities of Lemma 5.2, we have 
$$\left|a_i -1 \right| \  < \ {\epsilon \over n-r_i} \ \leq \ \epsilon \quad \text{for} \quad i=1, \ldots, m.$$

For a row in the $j$-th block of rows of type II the sum equals 
$$b_j=\sum_{i=1}^m {w_{ij} x_i y_j \over c_j  (1+w_{ij} x_i y_j)}.$$
By the inequalities of Lemma 5.2, we have 
$$|b_j -1 | \ < \  {\epsilon \over c_j} \ \leq \  \epsilon \quad \text{for} \quad j=1, \ldots, n.$$

Let us compute the column sums of $B(\epsilon)$. 
 
For the $j$-th column from the $i$-th block of columns the sum equals 
$$(n-r_i) {1 \over (n-r_i) (1+w_{ij} x_i y_j)} + c_j {w_{ij} x_i y_j \over c_j(1+w_{ij} x_i y_j)}=1.$$

Clearly, $B(\epsilon)$ is non-negative and hence by Lemma 4.3, we have 
$$\per B(\epsilon) \geq {(mn)! \over (mn)^{mn}} \phi(\epsilon) \quad \text{where} \quad 
\lim_{\epsilon \longrightarrow 0+} \phi(\epsilon)=1.$$
The proof now follows by (5.3.1) as $\epsilon \longrightarrow 0+$.
{\hfill \hfill \hfill} \qed

\head 6. Proofs of Lemmas 1.3 and 2.3 \endhead

We prove Lemma 2.3 only since Lemma 1.3 is a particular case of Lemma 2.3.

\demo{Proof of Lemma 2.3}
Since $H'(x)=\ln(1-x) -\ln x$, the value of the derivative at $x=0$ is $+\infty$ (we consider the 
right derivative there), the value of the derivative at $x=1$ is $-\infty$ (we consider 
the left derivative there) and the value of the derivative is finite for any $0< x< 1$.
Suppose that for the maximum entropy matrix $Z$ we have $z_{ij} \in \{0,1\}$ for some $i,j$ such that $w_{ij}=1$.
If $Y \in \PP(R, C; W)$, $Y=\left(y_{ij}\right)$,
 is a matrix such that $0 < y_{ij} <1$ whenever $w_{ij}=1$ then 
$$H\left( \epsilon Y + (1-\epsilon) Z \right) > H(Z) \quad 
\text{for a sufficiently small} \quad \epsilon >0,$$
which contradicts to the choice of $Z$. Hence 
$$0 < z_{ij} < 1 \quad \text{whenever} \quad w_{ij}=1.$$
Therefore, the gradient of $H(X)$ at $X=Z$ is orthogonal to the affine subspace of matrices 
$X=\left(x_{ij}\right)$ 
having row sums $R$, column sums $C$, and such that $x_{ij}=0$ whenever $w_{ij}=0$.
Hence
$$\ln {1-z_{ij} \over z_{ij}} =\lambda_i +\mu_j \quad \text{for all} \quad i,j \quad \text{such that}
\quad w_{ij}=1 \tag6.1$$ 
and some $\lambda_1, \ldots, \lambda_m$ and $\mu_1, \ldots, \mu_n$.
Hence
$$z_{ij}={e^{-\lambda_i} e^{-\mu_j} \over 1+e^{ -\lambda_i} e^{-\mu_j}} \quad \text{whenever} \quad w_{ij}=1.$$
Therefore,
$$\split &\sum \Sb j: \\ w_{ij}=1 \endSb {e^{-\lambda_i} e^{-\mu_j} \over 1+e^{ -\lambda_i} e^{-\mu_j}}  =
r_i \quad \text{for} \quad i=1, \ldots, m \\
&\sum \Sb i: \\ w_{ij}=1 \endSb {e^{-\lambda_i} e^{-\mu_j} \over 1+e^{ -\lambda_i} e^{-\mu_j}}  =
c_j \quad \text{for} \quad j=1, \ldots, n. \endsplit$$
Therefore, 
$$\sss^{\ast}=\left(-\lambda_1, \ldots, -\lambda_m \right) \quad 
\text{and} \quad \ttt^{\ast}=\left(-\mu_1, \ldots, -\mu_n \right)$$
is a critical point of 
$$G(\sss, \ttt; W) =-\sum_{i=1}^m r_i s_i -\sum_{j=1}^n c_j t_j +\sum \Sb (i,j): \\ w_{ij}=1 \endSb
\ln \left(1 + e^{s_i +t_j} \right).$$
Since $G$ is convex, $\left(\sss^{\ast}, \ttt^{\ast}\right)$ is also a minimum point.
Therefore, the point $\xx^{\ast}=\left(\xi_1, \ldots, \xi_m \right)$ and 
$\yy^{\ast}=\left(\eta_1, \ldots, \eta_n \right)$ where 
$$\xi_i=e^{-\lambda_i} \quad \text{for} \quad i=1, \ldots, m \quad \text{and} \quad
\eta_j=e^{-\mu_j} \quad \text{for} \quad j=1, \ldots, n$$
is a minimum point of 
$$F(\xx, \yy; W)=\left( \prod_{i=1}^m x_i^{-r_i} \right) \left(\prod_{j=1}^n y_j^{-c_j} \right)
\prod \Sb (i,j): \\ w_{ij}=1 \endSb \left(1+x_i y_j \right)$$
and satisfies 
$$\split &\sum \Sb j: \\ w_{ij}=1 \endSb {\xi_i \eta_j \over 1+\xi_i  \eta_j}  =
r_i \quad \text{for} \quad i=1, \ldots, m \\
&\sum \Sb i: \\ w_{ij}=1 \endSb {\xi_i \eta_j \over 1+\xi_i  \eta_j }  =
c_j \quad \text{for} \quad j=1, \ldots, n. \endsplit \tag6.2$$

Conversely, if $\xx^{\ast}=\left(\xi_1, \ldots, \xi_m \right)$ and $\yy^{\ast}=\left(\eta_1, \ldots, \eta_n\right)$
is a point where the minimum of $F(\xx, \yy; W)$ is attained, then, setting the gradient of 
$\ln F$ to 0, we obtain equations (6.2). Letting 
$$z_{ij} ={\xi_i \eta_j \over 1+ \xi_i \eta_j} \quad \text{when} \quad w_{ij} =1$$ 
and $z_{ij}=0$ when $w_{ij}=0$, we obtain a matrix $Z \in \PP(R, C; W)$. Moreover, the gradient of
$H(X)$ at $X=Z$ satisfies (6.1) with $\lambda_i =-\ln \xi_i$ and $\mu_j = -\ln \eta_j$, so 
$Z$ is the maximum entropy matrix. 

We now check:
$$\split H(Z)=&-\sum \Sb (i,j): \\ w_{ij} =1 \endSb z_{ij} \ln z_{ij} - \sum  \Sb (i,j): \\ w_{ij} =1 \endSb 
\left(1-z_{ij} \right) \ln \left(1 - z_{ij} \right) \\ = &
-\sum \Sb (i,j): \\ w_{ij} =1 \endSb {\xi_i \eta_j \over 1+\xi_i \eta_j} \ln {\xi_i \eta_j \over 1+\xi_i \eta_j} 
-\sum  \Sb (i,j): \\ w_{ij} =1 \endSb  {1 \over 1+\xi_i \eta_j} \ln {1 \over 1+ \xi_i \eta_j} \\
= &-\sum_{i=1}^m \ln \xi_i \left( \sum \Sb j: \\ w_{ij}=1 \endSb {\xi_i \eta_j \over 1+ \xi_i \eta_j} \right)
-\sum_{j=1}^n \ln \eta_j \left( \sum \Sb i: \\ w_{ij}=1 \endSb {\xi_i \eta_j \over 1+ \xi_i \eta_j} \right) \\
&\qquad\qquad + \sum \Sb (i,j): \\ w_{ij}=1 \endSb \ln \left(1 + \xi_i \eta_j \right) \\=
& -\sum_{i=1}^m r_i \ln \xi_i - \sum_{j=1}^n c_j \ln \eta_j +  \sum \Sb (i,j): \\ w_{ij}=1 \endSb \ln \left(1 + \xi_i \eta_j \right)
\endsplit$$
by (6.2). 
Hence 
$$H(Z) = \ln F\left(\xx^{\ast}, \yy^{\ast}; W \right)$$
and the proof follows.
\enddemo
{\hfill \hfill \hfill} \qed

\head 7. Proofs of Theorems 1.5 and 2.5 \endhead

We prove Theorem 2.5 only, since Theorem 1.5 is a particular case of Theorem 2.5. 
 
From formula (6.1), we have 
$${1 -z_{ij} \over z_{ij}} =e^{\lambda_i +\mu_j} \quad \text{for all} \quad i, j \quad \text{such that}
\quad w_{ij}=1$$ 
and some $\lambda_1, \ldots, \lambda_m$ and $\mu_1, \ldots, \mu_n$. Then, for all $i, j$ 
such that $w_{ij}=1$ and any 
$d_{ij} \in \{0, 1\}$, we have 
$$\split \Pr\bigl\{x_{ij}=d_{ij}\bigr\}=&z_{ij}^{d_{ij}} \left(1-z_{ij}\right)^{1-d_{ij}}=
\left(1 - z_{ij}\right) \left({1-z_{ij} \over z_{ij}} \right)^{-d_{ij}} \\ =
&\left(1-z_{ij} \right) e^{-(\lambda_i + \mu_j) d_{ij}}. \endsplit$$
Consequently, for any $D \in \Sigma(R, C; W)$, $D=\left(d_{ij}\right)$, we have 
$$\split \Pr\bigl\{X=D\bigr\}= &\prod_{i,j: \ w_{ij}=1} \left(1-z_{ij}\right)e^{-\left(\lambda_i + \mu_j\right)
d_{ij}} \\
=&\left(\prod_{i, j: \ w_{ij}=1} \left(1-z_{ij}\right)\right) \left( \prod_{i=1}^m e^{-\lambda_i r_i} \right) 
\left(\prod_{j=1}^n e^{-\mu_j c_j} \right). \endsplit$$
On the other hand,
$$\split e^{-H(Z)}=&\prod_{i, j: \ w_{ij}=1} z_{ij}^{z_{ij}} \left(1-z_{ij} \right)^{1-z_{ij}}\\=
&\left(\prod_{i, j: \ w_{ij}=1} \left(1-z_{ij}\right)\right) \left(\prod_{i, j: \ w_{ij} =1} \left({1-z_{ij} \over z_{ij}}
\right)^{-z_{ij}} \right) \\ =&\left(\prod_{i, j: \ w_{ij}=1} \left(1-z_{ij}\right)\right) 
\left( \prod_{i=1}^m e^{-\lambda_i r_i} \right) 
\left(\prod_{j=1}^n e^{-\mu_j c_j} \right), \endsplit$$
which completes the proof.
{\hfill \hfill \hfill} \qed

\head 8. Proofs of Theorems 1.4 and 2.4 \endhead

We prove Theorem 2.4 only since Theorem 1.4 is a particular case of Theorem 2.4.

We will use standard large deviation inequalities for bounded random variables, see, for example, 
Corollary 5.2 of \cite{Mc89}.

\proclaim{(8.1) Lemma} Let $Y_1, \ldots, Y_k$ be independent random variables such that 
$0 \leq Y_i \leq 1$ for $i=1, \ldots, k$. Let $Y=Y_1 + \ldots + Y_k$ and let $a=\EE Y$. 
Then, for $0 \leq \epsilon \leq 1$ we have 
$$\split &\Pr\bigl\{Y \ \geq \ (1+\epsilon) a \bigr\} \ \leq \ \exp\left\{-{1 \over 3} \epsilon^2 a \right\} 
\quad \text{and} \\ &\Pr\bigl\{ Y \ \leq \ (1-\epsilon) a \bigr\} \ \leq \ \exp\left\{ - {1 \over 2} \epsilon^2 a
\right\}. \endsplit$$
\endproclaim
{\hfill \hfill \hfill} \qed 

\subhead (8.2) Proof of Theorem 2.4 \endsubhead Let $X=\left(x_{ij}\right)$ be the $m \times n$
matrix of independent Bernoulli random variables such that $\EE X=Z$, as in Theorem 2.5.
By Theorem 2.5, the distribution of $X$ conditioned on $\Sigma(R, C; W)$ is uniform and hence 
$$\split \Pr \bigl\{D \in &\Sigma(R, C; W): \ \sigma_S(D) \leq (1-\epsilon) \sigma_S(Z)\bigr\} \\=
&{\Pr \bigl\{X: \ \sigma_S(X) \leq (1-\epsilon) \sigma_S(Z) \quad \text{and} \quad X \in \Sigma(R, C; W)
\bigr\} \over \Pr\bigl\{X: \ X \in \Sigma(R, C; W)\bigr\}}. \endsplit$$
Similarly,
$$\split \Pr \bigl\{D \in &\Sigma(R, C; W): \ \sigma_S(D) \geq (1+\epsilon) \sigma_S(Z)\bigr\} \\=
&{\Pr \bigl\{X: \ \sigma_S(X) \geq (1+\epsilon) \sigma_S(Z) \quad \text{and} \quad X \in \Sigma(R, C; W)
\bigr\} \over \Pr\bigl\{X: \ X \in \Sigma(R, C; W)\bigr\}}. \endsplit$$
By Theorem 2.5, Lemma 2.3 and Theorem 2.1, we get 
$$\split \Pr \bigl\{X \in \Sigma(R, C; W) \bigr\}\ =\ &e^{-H(Z)} |\Sigma(R, C; W)| \\ \geq \ 
&{(mn)! \over (mn)^{mn}} \left(\prod_{i=1}^m {(n-r_i)^{n-r_i} \over (n-r_i)!}\right) 
\left(\prod_{j=1}^n {c_j^{c_j} \over c_j!} \right) \\
\geq \ &(mn)^{-\gamma(m+n)} \endsplit$$
for some absolute constant $\gamma >0$.

Therefore,
$$\aligned \Pr \bigl\{D \in &\Sigma(R, C; W): \ \sigma(D) \leq (1-\epsilon) \sigma_S(Z) \bigr\} \\
\ \leq \ &(mn)^{\gamma(m+n)} \Pr \bigl\{X: \ \sigma_S(X) \leq (1-\epsilon) \sigma_S(Z)\bigr\} \\
\qquad \qquad &\text{and similarly} \\
 \Pr \bigl\{D \in &\Sigma(R, C; W): \ \sigma(D) \geq (1+\epsilon) \sigma_S(Z) \bigr\} \\
\ \leq \ &(mn)^{\gamma(m+n)} \Pr \bigl\{X: \ \sigma_S(X) \geq (1+\epsilon) \sigma_S(Z)\bigr\}.
\endaligned \tag8.2.1$$
By Lemma 8.1,
$$\aligned &\Pr\bigl\{X: \ \sigma_S(X) \ \leq \ (1-\epsilon)\sigma_S(Z) \bigr\} \ \leq \ \exp\left\{ -{1 \over 2} 
\epsilon^2 \sigma_S(Z) \right\} \\
&\qquad \qquad \text{and} \\ 
&\Pr\bigl\{X: \ \sigma_S(X) \ \geq \ (1+\epsilon)\sigma_S(Z) \bigr\} \ \leq \ \exp\left\{ -{1 \over 3} 
\epsilon^2 \sigma_S(Z) \right\}. \endaligned \tag8.2.2$$
Hence for 
$$\epsilon={\delta \ln n \over \sqrt{m}} \quad \text{and} \quad  \sigma_S(Z) \geq \delta mn$$
we have
$$\epsilon^2 \sigma_S(Z) \ \geq \ \delta^3  n \ln^2 n. \tag8.2.3$$ 
Combining (8.2.1)--(8.2.3), we conclude that for any $\kappa >0$ and all sufficiently large 
$n \geq m > q(\kappa, \delta)$ we have
$$\split &\Pr\bigl\{D \in \Sigma(R, C; W): \ \sigma_S(D) \leq (1-\epsilon) \sigma_S(Z) \bigr\} \ \leq \ 
n^{-\kappa n} \quad \text{and} \\
&\Pr\bigl\{D \in \Sigma(R, C; W): \ \sigma_S(D) \geq (1+\epsilon) \sigma_S(Z) \bigr\} \ \leq \ 
n^{-\kappa n} \endsplit$$
as required.
{\hfill \hfill \hfill} \qed

\head Acknowledgment \endhead

I am grateful to Alex Samorodnitsky for emphasizing the importance of matrix scaling in 
permanent computations during many enlightening conversations. I benefitted from conversations 
with Catherine Greenhill about asymptotic enumeration
 during the 2008 Schloss Dagstuhl meeting on design and analysis 
of randomized and approximation algorithms. After the first version of this paper was written,
John Hartigan pointed out to connections with the maximum entropy principle and suggested 
Theorem 1.5 to me (cf. \cite{BH09}), which led to a substantial simplification of the 
original proof and some strengthening of Theorems 1.4 and 2.4.


\Refs
\widestnumber\key{AAAAA}

\ref\key{Ba96}
\by A. Barvinok
\paper Two algorithmic results for the traveling salesman problem
\jour Mathematics of Operations Research
\vol 21 
\yr 1996
\pages  65--84
\endref

\ref\key{Ba07}
\by A. Barvinok
\paper Integration and optimization of multivariate polynomials by restriction onto a random subspace
\jour Foundations of Computational Mathematics 
\vol 7 
\yr 2007
\pages  229--244
\endref

\ref\key{Ba09}
\by A. Barvinok
\paper Asymptotic estimates for the number of contingency tables, integer flows, and volumes of transportation polytopes
\jour International Mathematics Research Notices. IMRN
\vol 2009
\yr 2009
\pages 348--385
\endref

\ref\key{BH09}
\by A. Barvinok and J. Hartigan
\paper Maximum entropy Gaussian approximation for the number of integer points and volumes of polytopes
\paperinfo preprint {\tt arXiv:0903.5223}
\yr 2009
\endref

\ref\key{Be74}
\by E.A. Bender
\paper The asymptotic number of non-negative integer matrices with given row and column sums
\jour Discrete Mathematics
\vol 10 
\yr 1974
\pages  217--223
\endref 

\ref\key{B+07}
\by I. Bez\'akov\'a, N. Bhatnagar, and E. Vigoda
\paper Sampling binary contingency tables with a greedy start
\jour  Random Structures $\&$ Algorithms 
\vol 30 
\yr 2007
\pages 168--205
\endref

\ref\key{C+08}
\by E.R. Canfield, C. Greenhill, and B.D. McKay
\paper Asymptotic enumeration of dense 0-1 matrices with specified line sums
\jour Journal of Combinatorial Theory. Series A 
\vol 115 
\yr 2008
\pages 32--66
\endref

\ref\key{CM05}
\by E. R. Canfield and B.D. McKay
\paper  Asymptotic enumeration of dense 0-1 matrices with equal row sums and equal column sums
\jour Electronic Journal of Combinatorics
\vol 12 
\yr 2005
\paperinfo Research Paper 29, 31 pp.
\endref

\ref\key{C+05}
\by Y. Chen, P. Diaconis, S.P. Holmes,  and J.S. Liu
\paper Sequential Monte Carlo methods for statistical analysis of tables
\jour  Journal of the American Statistical Association
\vol 100 
\yr 2005
\pages 109--120
\endref

\ref\key{CK09a}
\by W. Cuckler and J. Kahn
\paper Hamiltonian cycles in Dirac graphs 
\jour Combinatorica 
\vol 29 
\yr 2009
\pages 299--326
\endref 

\ref\key{CK09b}
\by W. Cuckler and J. Kahn
\paper Entropy bounds for perfect matchings and Hamiltonian cycles
\jour Combinatorica 
\vol 29 
\yr 2009
\pages  327--335
\endref

\ref\key{GC77}
\by I.J. Good and J.F. Crook
\paper The enumeration of arrays and a generalization related to contingency tables
\jour Discrete Mathematics
\vol 19 
\yr 1977
\pages 23--45
\endref

 \ref\key{G+06}
 \by C. Greenhill, B.D. McKay, and X. Wang
 \paper Asymptotic enumeration of sparse 0-1 matrices with irregular row and column sums
 \jour Journal of Combinatorial Theory. Series A 
 \vol 113 
 \yr 2006
 \pages 291--324
 \endref

\ref\key{GM09}
\by C. Greenhill and B.D.  McKay
\paper Random dense bipartite graphs and directed graphs with specified degrees
\jour Random Structures $\&$ Algorithms 
\vol 35 
\yr 2009
\pages 222--249
\endref

\ref\key{Gu04}
\by L. Gurvits
\paper Classical complexity and quantum entanglement
\jour Journal of Computer and System Sciences
\vol 69 
\yr 2004
\pages 448--484
\endref 

\ref\key{Gu08}
\by  L. Gurvits
\paper Van der Waerden/Schrijver-Valiant like conjectures and stable (aka hyperbolic) homogeneous polynomials: one theorem for all 
\jour Electronic Journal of Combinatorics
\paperinfo Research Paper 66
\vol 15(1)
\yr 2008
\endref

\ref\key{JS90}
\by M. Jerrum and A. Sinclair
\paper Fast uniform generation of regular graphs
\jour Theoretical Computer Science
\vol 73 
\yr 1990
\pages 91--100
\endref

\ref\key{J+04}
\by M. Jerrum, A. Sinclair, and E. Vigoda
\paper A polynomial-time approximation algorithm for the permanent of a matrix with nonnegative entries
\jour Journal of the ACM
\vol 51 
\yr 2004
\pages 671--697
\endref

\ref\key{Kh57}
\by A.I. Khinchin
\book Mathematical Foundations of Information Theory
\publ Dover Publications, Inc.
\publaddr New York, N. Y.
\yr 1957
\endref

\ref\key{L+00}
\by N. Linial, A. Samorodnitsky, and A.  Wigderson
\paper A deterministic strongly polynomial algorithm for matrix scaling and approximate permanents \jour Combinatorica 
\vol 20 
\yr 2000
\pages 545--568
\endref

\ref\key{LW01}
\by J.H. van Lint and R.M. Wilson
\book A Course in Combinatorics. Second edition
\publ Cambridge University Press
\publaddr Cambridge
\yr 2001
\endref

\ref\key{Ma95}
\by I.G. Macdonald
\book Symmetric Functions and Hall Polynomials. Second edition. With contributions by A. Zelevinsky
\bookinfo Oxford Mathematical Monographs
\publ  Oxford Science Publications. The Clarendon Press, Oxford University Press
\publaddr New York
\yr 1995
\endref

\ref\key{Mc89}
\by C. McDiarmid
\paper On the method of bounded differences
\inbook  Surveys in combinatorics, 1989 (Norwich, 1989)
\pages 148--188
\bookinfo London Mathematical Society Lecture Note Series
\vol 141
\publ Cambridge Univ. Press
\publaddr Cambridge
\yr 1989
\endref

\ref\key{Ne69}
\by P.E. O'Neil
\paper Asymptotics and random matrices with row-sum and column-sum restrictions
\jour Bulletin of the American Mathematical Society
\vol 75 
\yr 1969 
\pages 1276--1282
\endref

\ref\key{NN94}
\by Yu. Nesterov and A.  Nemirovskii
\book Interior-Point Polynomial Algorithms in Convex Programming
\bookinfo SIAM Studies in Applied Mathematics, 13
\publ  Society for Industrial and Applied Mathematics (SIAM)
\publaddr Philadelphia, PA
\yr 1994
\endref

\endRefs
\enddocument

\end